\documentclass[12pt,titlepage]{amsart}
\date{\today}
\usepackage{art}
\usepackage{amssymb}
\usepackage{eucal}

\newcommand{\Span}{\mathrm{Span}}
\newcommand{\mon}{\mathrm{Mon}}
\newcommand{\Gro}{Gr\"obner}

\newcommand{\bld}{block decomposition}
\newcommand{\blds}{block decompositions}

\begin{document}
\bibliographystyle{plain}

\begin{titlepage}
\vspace{2 in}
\begin{center} BLOCK STANLEY DECOMPOSITIONS\\
I. ELEMENTARY AND GNOMON DECOMPOSITIONS  \end{center}
\vspace{0.5in}

\begin{center}
James Murdock\\
Department of Mathematics\\
Iowa State University\\
Ames, Iowa 50011\\
jmurdock@iastate.edu\\
\vspace{.2in}
and\\
\vspace{.2in}
Theodore Murdock\\
823 Carroll Ave.\\
Ames, Iowa 50010\\
theodoremurdock@gmail.com
\end{center}
\vspace{1in}
This is the final preprint.  The published version appears in J. of Pure and Applied Algebra 219 (2015) 2189-2205.  The doi is 10.1016/j.jpaa.2014.07.030.

\vspace{1in} 
Running head: BLOCK STANLEY DECOMPOSITIONS I
\vspace{1in}
\begin{center}
Corresponding Author\\
James Murdock\\
phone: 515-232-7945\\
email: jmurdock@iastate.edu\\
fax: 515-294-5454
\end{center}
\end{titlepage}

ABSTRACT:  Stanley decompositions are used in invariant theory and the theory of normal forms for dynamical systems to provide a unique way of writing each invariant as a polynomial in the Hilbert basis elements.  Since the required Stanley decompositions can be very long, we introduce a more concise notation called a block decomposition, along with three notions of shortness (incompressibility, minimality of Stanley spaces, and minimality of blocks) for block decompositions.  We give two algorithms that generate different block decompositions, which we call elementary and gnomon decompositions, and give examples.  Soleyman-Jahan's criterion for a Stanley decomposition to come from a prime filtration is reformulated to apply to block decompositions.  We simplify his proof, and apply the theorem to show that elementary and gnomon decompositions come from ``subprime'' filtrations.  In a sequel to this paper we will introduce two additional algorithms that generate block decompositions that may not always be subprime, but are always incompressible.

\newpage
\markright{BLOCK STANLEY DECOMPOSITIONS I}

\section{\bf Introduction}\label{intro}

Let $I\subset S=K[\bx]=K[x_1,\dots,x_n]$ be a monomial ideal and let $M$ be the set of {\it standard monomials}, that is, the monomials $\bx^{\bk}=x_1^{k_1}\cdots x_n^{k_n}\in K[\bx]$ that are {\it not} in $I$.  Each coset in $S/I$ has exactly one representative as a {\it standard polynomial} in $\Span(M)$, and $S/I$ is isomorphic to $\Span(M)$ as a vector space over $K$, although the ring structure of $S/I$ is lost.  A {\it \bld} of $M$  (defined precisely in Section \ref{block}) is an expression for $M$ as a disjoint union of rectangular {\it blocks} of standard monomials, regarded as points in Newton space.  A block will be called a {\it Stanley block} if its span is a Stanley space; a \bld\ of $A$ into Stanley blocks is equivalent to a \sd\ of $S/I$.

In this paper and its sequel we study algorithms leading from $I$ to a \bld\ for $M$, emphasizing that different algorithms produce \blds\ with different algebraic and geometric properties.  Our long-range goal is to produce \blds\ that are as simple as possible, with respect to various criteria.  In this first paper we present two algorithms, producing what we call the {\it elementary} and {\it gnomon} \blds\ respectively, and show that these are associated with {\it subprime filtrations} of $I$, a natural generalization of prime filtrations.  In doing so, we generalize Soleyman-Jahan's criterion in \cite[Prop. 2.7]{AliPrime} for a \sd\ to be associated with a prime filtration, and simplify the proof.  The gnomon decomposition is shorter (contains fewer blocks) than the elementary decomposition.  In the second paper we present algorithms leading to what we call the {\it organized} and {\it stacked} \blds.   These are {\it incompressible}, in the sense that they cannot be shortened by combining blocks.  We have not yet found an algorithm that always produces a minimal \bld\ (one with the fewest blocks), but since any such \bld\ must be incompressible, these algorithms are a step in that direction.

From one point of view, \blds\ are a generalization of \sds, and are in general {\it coarser} than \sds.  From another point of view, each \bld\ has a unique minimal Stanley refinement, given by Algorithm \ref{SDfromBD} below, and since this algorithm is very easy, the \bld\ can be taken as a {\it compressed notation} for this \sd.  From this viewpoint, we are studying algorithms leading from $I$ to a \sd\ of $S/I$.  The usual algorithm of this type, \cite[Lemma 2.4]{SturmWh}, proceeds by induction on the variables $x_1,\dots,x_n$.  Our gnomon algorithm works by induction on the generators of $I$, so that if one has determined a \bld\ for $S/I$ and then adds a generator to $I$, it is not necessary to start over from the beginning.

In Section \ref{block} we define \blds\ and \sds\ and present Algorithm \ref{SDfromBD} that relates them. Sections \ref{elem} and \ref{gnomon} contain the algorithms for the elementary and gnomon block decompositions.  Section \ref{prime} contains the generalization of Soleyman-Jahan's theorem and the proof that elementary and gnomon decompositions are subprime.

Stanley decompositions are best known among algebraists in connection with the conjecture that the Stanley depth of a module over a polynomial algebra is an upper bound for the classical depth, but have a separate life in applied mathematics through their use in describing subalgebras.  This life is based entirely on the first two pages of \cite{SturmWh} and, for one  application, on \cite{SturmWhBracket}.  The motivation for our work comes from that direction.  Let $\CA$ be a subalgebra of $K[\by]=K[y_1,\dots,y_m]$.  It is desired to find a formula, containing finitely many arbitrary polynomial functions, that expresses each element of $\CA$ exactly once as the arbitrary functions are varied.  A set of polynomials $f_1,\dots,f_n\in\CA$ is a {\it \hb} for $\CA$ if the monomials $f_1^{k_1}\cdots f_n^{k_n}$  span $\CA$ as a vector space over $K$.  The expression for a given $f$ in terms of the \hb\ is usually not unique, because the monomials $f_1^{k_1}\cdots f_n^{k_n}$ may not be linearly independent.  There always exist linearly independent subsets of these monomials, which we call {\it preferred sets}, and the expression of any $f$ as a linear combination of these is unique.  But a preferred set of monomials is infinite, and may not have a finite description leading to an expression involving arbitrary polynomial functions.  Stanley decompositions in the manner of \cite{SturmWh} provide suitable preferred sets . Let $\Phi:K[\bx]\ra K[\by]$ be the unique algebra homomorphism such that $\Phi(x_i)=f_i$; then $J=\ker\Phi$ is the ideal of relations among the $f_i$.  Compute a \Gro\ basis for $J$ using an appropriate elimination order on $K[\bx,\by]$ (\cite[p. 81]{AdLou}), and obtain the monomial ideal $I$ of leading terms of $J$, which is generated by the leading terms of the \Gro\ basis.  The standard monomials of $I$ provide a preferred set of monomials in $f_1,\dots,f_n$ when each $x_i$ is replaced by $f_i$.  In fact $\CA$ is isomorphic as an algebra to $K[\bx]/J$, and isomorphic as a vector space to both $K[\bx]/I$ and $\Span(M)$.  Furthermore, a \sd\ for $K[\bx]/I$ provides the desired formula expressing each element of $\CA$ uniquely.  For instance, if $\CA$ has a \hb\ of four elements $f_1,\dots,f_4$ with the single relation $f_3^2-f_2^3-f_1^2f_4=0$, and $f_3^2$ is taken as the leading term, then $I=\langle x_3^2\rangle$ and the \sd\ given by the algorithm of \cite{SturmWh} is
\begin{equation}\label{N4}
K[x_1,x_2,x_3,x_4]/\langle x_2^3\rangle \cong K[x_1,x_2,x_4]\oplus K[x_1,x_2,x_4]x_3.
\end{equation}
It follows that every element of $\CA$ can be written uniquely in the form
\[
f=F_1(f_1,f_2,f_4)+F_2(f_1,f_2,f_4)f_3,
\]
where $F_1$ and $F_2$ are arbitrary polynomials in three variables.  (This example arises for the  classical seminvariants of the binary cubic form.)

Notice that the number of arbitrary polynomials is the same as the number of Stanley spaces.  The number of such polynomials can become large in larger problems, and it is natural to ask if this can be reduced.  In part this can be done by a good choice of the term order; for instance, if $f_2^3$ was taken as the leading term in the example, there would be three Stanley spaces instead of two. But even once $I$ is selected, there can be large differences in the number of Stanley terms in different \sds\ for $S/I$.  This suggests the problem:  {\it Find a minimal \sd\ for $S/I$ given a fixed monomial ideal $I$.}   A related, but different, problem is to {\it find a minimal \bld\ for $M$ given $I$}; here the minimality is with regard to the number of blocks, not the number of Stanley spaces.

It is not an accident that the example mentioned above comes from classical invariant theory.  Invariant theory is closely linked with normal form theory for systems of differential equations with nilpotent linear part.  There was a flurry of work on this topic in the late 1980s and early 1990s (\cite{CSNilp}, \cite{CSN4}, \cite{CSW}, \cite{SturmWhBracket}, \cite{CSSurv}, \cite{SturmWh}), and a second flurry more recently (\cite{MurdStruc}, \cite{MurdNTA}, \cite[Ch.12]{SanVerMurd}, \cite{SanJoint}, \cite{Mal}, \cite{Mal-2}, \cite{GachInv}, \cite{GachNF}).  Additional motivation for our questions comes from these papers.  The algebra of scalar fields that are invariant under the the group $\{e^{Mt}:t\in\R\}$, where $M$ is a nilpotent matrix, plays an important role.  These invariants, in turn, are the same as classical seminvariants of a binary form, or joint seminvariants of several forms.  Knowing the seminvariants of several forms, {\it transvectants} are used to find a \hb\ for the the joint seminvariants.  The {\it box product} method introduced in \cite{MurdNTA} lifts this procedure to the level of \sds.  The box product of two \sds\ for the seminvariants of two binary forms is a \sd\ for the joint seminvariants of the two forms; thus the uniqueness issue is taken care of automatically, without further \Gro\ basis work.  In \cite{MurdBox} I (J.M.) have reformulated the box product as a product of block decompositions, in which $\boxtimes$ distributes over disjoint union; the basic unit of computation is the box product of two blocks, for which algorithms are given, and the final result (expressed in \sds) is shorter than for the method in \cite{MurdNTA}.  It is best, then, to start with \blds\ with a minimal number of blocks.

This research began when I (J.M.) asked T.M., a software engineer, to write programs for me concerning \sds.  He created block decompositions as a computer-friendly shorthand for \sds, and the gnomon algorithm to provide an interactive program in which ideal generators could be added one at a time.  I provided the mathematical write-up, proofs, and examples, and the material on subprime filtrations.  (The word {\it subprime} was provided by the mortgage crisis.)

\section{\bf Block decompositions}\label{block}
Let $V\subset S=K[\bx]$ be a vector subspace; $V$ is a {\bf monomial space} if it is spanned by the set $M=\mon(V)$ of monomials in $V$, which is then also the unique monomial basis for $V$.   Under the one-to-one correspondence $V=\Span(M)$, $M=\mon(V)$ between monomial spaces and their monomial bases, a direct sum of two spaces corresponds to the disjoint union of their bases.  We often identify $\bx^{\bbm}$ with $\bbm=(m_1,\dots,m_n)$ in the {\bf Newton space} $\N^n$.

A {\bf Stanley space} is a monomial space of the form
\[
K[X]\p = \{f\p:f\in K[X]\},
\]
where $X$ is a subset of the set of variables $\{x_1,\dots,x_n\}$ and $\p$ is a monomial.  An expression exhibiting a monomial space $V$ as a direct sum of finitely many Stanley spaces is called a {\bf Stanley decomposition} of $V$.  (Because of this finiteness requirement, there are monomial spaces, such as $\Span\{1,x^2,x^4,x^6,\dots\}\subset K[x]$, that do not have Stanley decompositions.)

Let $I$ be a monomial ideal in $S$ and let $M=\mon(S)\smallsetminus \mon(I)$ be the set of {\bf standard monomials} for $I$.  Then $V=\Span(M)$ is the space of {\bf standard polynomials}, and as noted in the introduction, $S/I\cong V$ as vector spaces.  If $S=K[x,y]$ and $I=\langle xy^3, x^3y\rangle$, a Stanley decomposition for $S/I$ is
\begin{equation}\label{sd1}
K[x]\oplus K[y]y \oplus Kxy\oplus Kx^2y\oplus Kxy^2\oplus Kx^2y^2.
\end{equation}
Applying the correspondences above, we replace each Stanley space by its set of monomials and replace $\oplus$ by the disjoint union symbol $\sqcup$, obtaining
\[
\{1,x,x^2,\dots\}\sqcup\{y,y^2,y^3,\dots\}\sqcup\{xy, xy^2, x^2y, x^2y^2\}.
\]
In this way the \sd\ (\ref{sd1}) becomes the following {\it block decomposition} in Newton space:

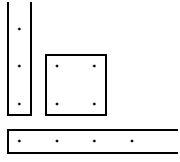
\begin{center}
\setlength{\unitlength}{0.5mm}
\begin{picture}(100,30)(-40,0)
\put(0,0){\circle*{1}}
\put(10,0){\circle*{1}}
\put(20,0){\circle*{1}}
\put(30,0){\circle*{1}}
\put(0,10){\circle*{1}}
\put(0,20){\circle*{1}}
\put(0,30){\circle*{1}}
\put(10,10){\circle*{1}}
\put(10,20){\circle*{1}}
\put(20,10){\circle*{1}}
\put(20,20){\circle*{1}}
\put(7,7){\line(1,0){16}}
\put(7,23){\line(1,0){16}}
\put(7,7){\line(0,1){16}}
\put(23,7){\line(0,1){16}}
\put(-3,-3){\line(1,0){45}}
\put(-3,-3){\line(0,1){6}}
\put(-3,3){\line(1,0){45}}
\put(-3,7){\line(1,0){6}}
\put(-3,7){\line(0,1){30}}
\put(3,7){\line(0,1){30}}
\end{picture}
\end{center}

More generally, let $[a,b]$ denote an integer interval, with the understanding that $[a,\infty]$ does not include $\infty$ and that $[a,b]=\emptyset$ if $b<a$.  Define a {\bf block} $B$ in $\N^n$ to be a Cartesian product of $n$ intervals:
\begin{equation}\label{defblock}
B=\bm \bb\\ \ba\ebm = \bm b_1&b_2&\cdots&b_n\\a_1&a_2&\cdots&a_n\ebm = [a_1,b_1]\times\cdots\times [a_n,b_n].
\end{equation}
If any $b_i<a_i$, the block is empty.  The bottom row $\ba$ of a block is called its {\bf inner corner}.  The letter $B$ will be used ambiguously to mean the matrix as such, the set of integer vectors represented by the matrix, and the set of monomials having those integer vectors as exponents.  Let $M\subset\mon(S)$.  A {\bf block decomposition} of $M$ is an expression
\begin{equation}\label{bd}
M=B^1\sqcup B^2\sqcup \dots \sqcup B^s,
\end{equation}
where each $B^k$ is a block.
For example, the Stanley decomposition (\ref{sd1}) can be written as the block decomposition
\begin{equation}\label{bd1}
B^1\sqcup B^2\sqcup B^3 = \bm \infty&0\\0&0\ebm \sqcup \bm 0&\infty\\0&1\ebm \sqcup \bm 2&2\\1&1\ebm,
\end{equation}
corresponding to the diagram shown above in Newton space.

\brem After this research was mostly completed, we learned of the {\bf interval decompositions} defined in \cite{HerVladZh}.  They define the {\bf interval} $[\ba,\bb]$ to be the set of monomials lying between $\bx^{\ba}$ and $\bx^{\bb}$ in division partial order.  This is essentially the same as our block $B$ in (\ref{bd1}), except that they consider only monomials that divide a chosen $\bx^{\bc}$, so that, in particular, $\infty$ is not allowed in the top row.  This is sufficient since they only consider interval decompositions for $I/J$, where $I$ and $J$ are two monomial ideals in $S$.\erem

A block decomposition is called {\bf compressible} if the union of some subset of the blocks is itself a block.  In this case the decomposition can be simplified by performing the union.  However, performing one such union may prevent another one, so a given decomposition may be compressible in several ways to give distinct incompressible decompositions.  A decomposition is called {\bf minimal} if there is no decomposition (describing the same set of monomials) having fewer blocks.  Incompressibility and minimality are two possible notions of ``simplicity'' of a decomposition (as briefly discussed in the introduction).  An incompressible decomposition need not be minimal, and a minimal decomposition need not be unique.  The following example shows that incompressibility is a global property of a decomposition, and cannot be detected by examining the decomposition locally:  The disjoint union
\[
\bm 1&0\\0&0\ebm\sqcup\bm 0&2\\0&1\ebm\sqcup\bm1&1\\1&1\ebm\sqcup\bm2&1\\2&0\ebm\sqcup\bm2&2\\1&2\ebm
\]
is compressible to the single block $\bm 2&2\\0&0\ebm$, although no proper subset of the five blocks can be compressed.

A {\bf Stanley block} is a block (\ref{defblock}) in which each $b_i$ equals either $a_i$ or $\infty$.  The span of a Stanley block is a Stanley space, and any Stanley decomposition can be changed to a block decomposition by replacing Stanley spaces by Stanley blocks and $\oplus$ by $\sqcup$.  Such a block decomposition can almost always be compressed.  For example, $K[x_1,x_2,x_4]\oplus K[x_1,x_2,x_4]x_3$ becomes
\[
\bm \infty&\infty&0&\infty\\ 0&0&0&0\ebm \sqcup \bm \infty&\infty&1&\infty\\0&0&1&0\ebm=\bm \infty&\infty&1&\infty\\0&0&0&0\ebm.
\]

Any block $B$ can be written as a disjoint union of its {\bf bounded part} and its {\bf unbounded part}; the bounded part is obtained by replacing the unbounded columns with zeroes, the unbounded part by doing the same with the bounded columns.  For instance,
\[
B=\bm \infty&5&\infty&2\\1&3&0&1\ebm = \bm 0&5&0&2\\0&3&0&1\ebm \sqcup \bm \infty&0&\infty&0\\1&0&0&0\ebm.
\]
Let $X$ be the set of variables associated with unbounded columns of $B$ ($x_1$ and $x_3$ in the example).  Let $\psi$ be the monomial in the bottom row of the unbounded part ($x_2^2x_4$ in the example).  Let $\th_1,\dots,\th_s$ be the finite set of monomials in the bounded part ($x_2^3x_4$, $x_2^4x_4$, $x_2^5x_4$, $x_2^3x_4^2$, $x_2^4x_4^2$, $x_2^5x_4^2$, with $s=6$).  Let $\p_k=\th_k\psi$ for $k=1,\dots,s$.  Then the block $B$ is converted to its minimal \sd\ as follows.
\balg\label{SDfromBD}   The Stanley decomposition with the fewest Stanley spaces that describes the span of the monomials in $B$ is
\[
\Span(B)=K(X)\p_1\oplus\cdots\oplus K(X)\p_s.
\]
\ealg

\noindent Note that all of the coefficient algebras $K(X)$ in the minimal \sd\ of a block are equal, and that the number of Stanley spaces $s$ is the product of the numbers of elements in the intervals defined by the bounded columns of $B$.  When Algorithm \ref{SDfromBD} is applied to each block in a \bld, the coefficient rings may differ from block to block, but the total number of Stanley spaces is just the sum of the number for each block.

As an introduction to the next two sections, we now give examples of the {\it elementary} and {\it gnomon} block decompositions defined in those sections.  For
\begin{equation}\label{firstexample0}
I=\langle x^3y^9, x^7y^5\rangle\subset S=K[x,y],
\end{equation}
the {\it elementary} decomposition of $S/I$ is the disjoint union of the blocks
\begin{equation}\label{firstexample1}
\bm 2&4\\0&0\ebm, \bm 2&8\\0&5\ebm, \bm 2&\infty\\0&9\ebm,\bm 6&4\\3&0\ebm, \bm 6&8\\3&5\ebm, \bm \infty&4\\7&0\ebm.
\end{equation}
The {\it gnomon} decomposition comes in two forms,
\begin{equation}\label{firstexample2}
\bm 2&\infty\\0&0\ebm\sqcup \bm 6&8\\3&0\ebm\sqcup \bm \infty&4\\7&0\ebm
\end{equation}
and
\begin{equation}\label{firstexample3}
\bm \infty&4\\0&0\ebm\sqcup \bm6&\infty\\0&5\ebm\sqcup \bm 2&\infty\\0&9\ebm.
\end{equation}
The elementary decomposition is compressible to either of the gnomon decompositions, and (in this case, but not always) the gnomon decompositions are incompressible.

\section{The elementary block decomposition}\label{elem}

In Newton space, the {\bf division partial order} $\bx^{\bbm}\vert \bx^{\bbm'}$ will be written as $\bbm \preceq \bbm'$, meaning $m_i$ $\leq m'_i$ for $i=1,\dots,n$.  For variables in Newton space we use $\bmu=(\m_1,\dots,\m_n)$, reserving  $\bbm=(m_1,\dots,m_n)$ for constants, so that an equation such as $\m_2=m_2$ will represent the ``hyperplane'' through $\bbm$ perpendicular to the $\m_2$ axis in $\N^n$.  Let $S=K[\bx]$ and let $I=\langle \bx^{\bbm^1},\dots,\bx^{\bbm^r}\rangle \subset S$ be the monomial ideal with the indicated generators; we also write $I=\langle \bbm^1,\dots,\bbm^r\rangle$.  It is assumed that these generators are minimal, that is, no redundant generators (divisible by another generator) are included.  The {\bf elementary block decomposition} for $S/I$ is created by first gridding the Newton space $\N^n$ with the hyperplanes $\m_i=m^r_i$ for $i=1,\dots,n$ and $j=1,\dots,s$, passing through all points of the minimal generating set, and then discarding those of the resulting blocks that belong to $\mon(I)$ rather than $\mon(S)\smallsetminus \mon(I)$.  The following figure illustrates the elementary decomposition (\ref{firstexample1}) at the end of Section \ref{block}.  The solid dots are the generators at $(3,9)$ and $(7,5)$, the dark lines are the boundary of $\mon(I)$, and the blocks are numbered in the order that they appear in (\ref{firstexample1}).

\begin{center}
\setlength{\unitlength}{2mm}
\begin{picture}(50,15)(-18,0)


\put(0,0){\line(1,0){12}}
\put(0,0){\line(0,1){15}}


\put(1,-0.25){\line(0,1){0.5}}
\put(2,-0.25){\line(0,1){0.5}}
\put(3,-0.25){\line(0,1){0.5}}
\put(4,-0.25){\line(0,1){0.5}}
\put(5,-0.25){\line(0,1){0.5}}
\put(6,-0.25){\line(0,1){0.5}}
\put(7,-0.25){\line(0,1){0.5}}
\put(8,-0.25){\line(0,1){0.5}}
\put(9,-0.25){\line(0,1){0.5}}
\put(10,-0.25){\line(0,1){0.5}}
\put(11,-0.25){\line(0,1){0.5}}

\put(-0.25,1){\line(1,0){0.5}}
\put(-0.25,2){\line(1,0){0.5}}
\put(-0.25,3){\line(1,0){0.5}}
\put(-0.25,4){\line(1,0){0.5}}
\put(-0.25,5){\line(1,0){0.5}}
\put(-0.25,6){\line(1,0){0.5}}
\put(-0.25,7){\line(1,0){0.5}}
\put(-0.25,8){\line(1,0){0.5}}
\put(-0.25,9){\line(1,0){0.5}}
\put(-0.25,10){\line(1,0){0.5}}
\put(-0.25,11){\line(1,0){0.5}}
\put(-0.25,12){\line(1,0){0.5}}
\put(-0.25,13){\line(1,0){0.5}}
\put(-0.25,14){\line(1,0){0.5}}


\put(3,9){\circle*{1}}
\put(7,5){\circle*{1}}
\thicklines
\put(3,9){\line(0,1){6}}
\put(3,9){\line(1,0){4}}
\put(7,5){\line(0,1){4}}
\put(7,5){\line(1,0){5}}


\thinlines
\put(2,0){\line(0,1){4}}
\put(0,4){\line(1,0){2}}

\put(0,5){\line(1,0){2}}
\put(2,5){\line(0,1){3}}
\put(0,8){\line(1,0){2}}

\put(0,9){\line(1,0){2}}
\put(2,9){\line(0,1){6}}

\put(3,0){\line(0,1){4}}
\put(6,0){\line(0,1){4}}
\put(3,4){\line(1,0){3}}

\put(3,5){\line(0,1){3}}
\put(6,5){\line(0,1){3}}
\put(3,5){\line(1,0){3}}
\put(3,8){\line(1,0){3}}

\put(7,0){\line(0,1){4}}
\put(7,4){\line(1,0){5}}


\put(1,1.5){1}
\put(1,5.5){2}
\put(1,11){3}
\put(4,1.5){4}
\put(4,5.5){5}
\put(9,1.5){6}
\end{picture}
\end{center}

The elementary block decomposition is almost always compressible, and so is not very desirable in itself, but it shows that a natural block decomposition always exists.  The other decompositions (gnomon, organized, and stacked) in this paper and its sequel can be obtained by compressing the elementary decomposition, so we will never need to consider blocks smaller than the blocks in the elementary decomposition, or blocks that are not disjoint unions of elementary blocks.

The following is a precise algorithm to create the elementary decomposition.

\balg\label{ebd} Let $I=\langle \bbm^1,\dots,\bbm^r\rangle$ be presented by its minimal generating set.  The elementary block decomposition is generated as follows.
\blist
\item For $i=1,\dots,n$, create the list of exponents of $x_i$ in the set of generators, adding zero at the beginning of each list:
\begin{align*}
&L_1:  0, m^1_1,\dots, m^r_1;\\
&L_2:  0, m^1_2,\dots, m^r_2;\\
&\vdots\\
&L_n:  0, m^1_n,\dots,m^r_n.
\end{align*}
\item\label{reorder} Delete any repetitions in each list $L_i$, and re-order each list in increasing numerical order.
\item Create a preliminary list of inner corners as follows.  For $i=1,\dots,n$, choose an entry $a_i$ from $L_i$; then create an inner corner $\ba=(a_1,\dots,a_n)$.  Do this in all possible ways.
\item\label{refine} Refine the list of inner corners by discarding any $\ba\in I$, that is, any $\ba$ such that $\bbm^i\prec \ba$ for some $i$.
\item\label{createblock} Create a block $B$ for each inner corner $\ba$ in the refined list, by determining its outer corner $\bb$, as follows.  If $a_i$ has a successor $a_i'$ in $L_i$ (under the ordering from step \ref{reorder}), put $b_i=a_i'-1$.  If $a_i$ is the last element of $L_i$, but $b_i=\infty$.
\item\label{elemenum} Let $s$ be the number of resulting blocks, and enumerate these as
\[
B^k=\bm \bb^k\\ \ba^k \ebm
\]
for $k=1,\dots,s$.
\elist
\ealg

\brem\label{elemorder} The order of enumeration in Step \ref{elemenum} does not matter except in Section \ref{prime}.  For that section, enumerate the blocks in lexicographic order by their inner corners.\erem

For the example (\ref{firstexample0}), the lists from Step 2 are $L_1=\{0,3,7\}$, $L_2=\{0,5,9\}$.  The preliminary list of inner corners (in lexicographic order) is $\{(0,0), (0,5),(0,9),(3,0),(3,5)(3,9),(7,0),(7,5),(7,9)\}$.  From this we discard $(3,9)$, $(7,5)$, and $(7,9)$.  It is now easy to check that the outer corners in (\ref{firstexample1}) follow from step 5 above.

\brem\label{elemprop} The elementary decomposition has the following property, which will be used in Section \ref{prime}.  Let $\bc \notin I$ be a monomial.  Then there exists a unique inner corner $\ba^k$ of the elementary decomposition for $S/I$ that is maximal (under $\prec$) among all inner corners satisfying $\ba^k\preceq\bc$.  The inner corner $\ba^k$ satisfying this condition determines the block $B^k$ that contains $\bc$.  (The proof is trivial: The block containing $\bc$ has for the $i$-th component of its inner corner the largest element of $L_i$ that is $\leq c_i$.)\erem

\section{The gnomon decomposition}\label{gnomon}

Consider the following {\bf Block Subtraction Problem:}   Given a block (\ref{defblock}) and a principal monomial ideal $I=\langle \bbm \rangle$, find a block decomposition of $B\smallsetminus I$.  If $\bbm$ is within the block $B$, there are $n!$ distinct natural solutions.  For $n=2$, the following figure shows a block $B$ in $\N^2$, an ideal $I$ generated by the heavy dot, and the two block decompositions of $B\smallsetminus I$.  The set of points in $B\smallsetminus I$ is the sort of L-shaped region the Pythagoreans called a  {\bf gnomon} ($\g\n\o\m\o\n$).

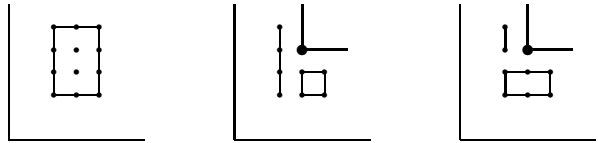
\begin{center}
\setlength{\unitlength}{3mm}
\begin{picture}(50,6)(-5,0)


\put(5,0){\line(1,0){6}}
\put(5,0){\line(0,1){6}}

\put(7,2){\circle*{0.25}}
\put(7,3){\circle*{0.25}}
\put(7,4){\circle*{0.25}}
\put(7,5){\circle*{0.25}}

\put(8,2){\circle*{0.25}}
\put(8,3){\circle*{0.25}}
\put(8,4){\circle*{0.25}}
\put(8,5){\circle*{0.25}}

\put(9,2){\circle*{0.25}}
\put(9,3){\circle*{0.25}}
\put(9,4){\circle*{0.25}}
\put(9,5){\circle*{0.25}}

\thinlines
\put(7,2){\line(1,0){2}}
\put(7,2){\line(0,1){3}}
\put(7,5){\line(1,0){2}}
\put(9,2){\line(0,1){3}}


\put(15,0){\line(1,0){6}}
\put(15,0){\line(0,1){6}}

\put(17,2){\circle*{0.25}}
\put(17,3){\circle*{0.25}}
\put(17,4){\circle*{0.25}}
\put(17,5){\circle*{0.25}}

\put(18,2){\circle*{0.25}}
\put(18,3){\circle*{0.25}}
\put(18,4){\circle*{0.5}}

\put(19,2){\circle*{0.25}}
\put(19,3){\circle*{0.25}}

\thicklines
\put(18,4){\line(1,0){2}}
\put(18,4){\line(0,1){2}}

\thinlines
\put(17,2){\line(0,1){3}}
\put(18,2){\line(1,0){1}}
\put(18,3){\line(1,0){1}}
\put(19,2){\line(0,1){1}}
\put(18,2){\line(0,1){1}}


\put(25,0){\line(1,0){6}}
\put(25,0){\line(0,1){6}}

\put(27,2){\circle*{0.25}}
\put(27,3){\circle*{0.25}}
\put(27,4){\circle*{0.25}}
\put(27,5){\circle*{0.25}}

\put(28,2){\circle*{0.25}}
\put(28,3){\circle*{0.25}}
\put(28,4){\circle*{0.5}}

\put(29,2){\circle*{0.25}}
\put(29,3){\circle*{0.25}}

\thicklines
\put(28,4){\line(1,0){2}}
\put(28,4){\line(0,1){2}}

\thinlines
\put(27,2){\line(0,1){1}}
\put(27,2){\line(1,0){2}}
\put(27,3){\line(1,0){2}}
\put(29,2){\line(0,1){1}}
\put(27,4){\line(0,1){1}}

\end{picture}
\end{center}

\noindent The first of these decompositions, expressed algebraically, takes the form
\begin{equation}\label{leftsolution}
\bm b_1&b_2\\a_1&a_2\ebm \smallsetminus \langle (m_1,m_2) \rangle = \bm m_1-1&b_2\\a_1&a_2\ebm \sqcup \bm b_1&m_2-1\\m_1&a_2\ebm.
\end{equation}
We will generalize this formula to arbitrary $n$ in such a way that it gives all $n!$ solutions, one for each permutation of the variables $x_1,\dots,x_n$.  First we focus on the natural order of variables (the trivial permutation), which for $n=2$ gives (\ref{leftsolution}).

The following solution for $n=2$ is valid whether or not $\bbm$ is in the interior of $B$, and reduces to (\ref{leftsolution}) when it is:
\begin{align*}
&\bm b_1&b_2\\a_1&a_2\ebm \smallsetminus \langle (m_1,m_2)\rangle =\\
&\qquad \bm \min\{b_1,m_1-1\}&b_2\\a_1&a_2\ebm \sqcup \bm b_1&\min\{b_2,m_2-1\}\\\max\{a_1,m_1\}&a_2\ebm.
\end{align*}
The first block consists of all monomials in $B$ that are not divisible by $x_1^{m_1}$, and is empty if $m_1\leq a_1$.  The second block consists of those monomials in $B$ that {\it are} divisible by $x_1^{m_1}$ but not by $x_2^{m_2}$.  But this solution has a defect: it sometimes divides $B$ unnecessarily into two nonempty blocks whose union is again $B$.  This happens whenever $a_1<m_1<b_1$ but $b_2<m_2$.  In fact, if $m_i>b_i$ for $i=$ 1 or 2 or both, then $B\cap I =\emptyset$ and $B\smallsetminus I=B$.  This observation also simplifies the formula in the remaining cases (where $m_i\leq b_i$ for all $i$), because then $\min\{b_i,m_i-1\}$ can be replaced by $m_i-1$.

Generalizing to higher dimensions, we define the {\bf gnomon} determined by $B$ and $I$ to be the set represented by the $3\times n$ matrix
\[
\bm \bb \\ \bbm \\ \ba \ebm = \bm b_1&\cdots&b_n\\m_1&\cdots&m_n\\a_1&\cdots&a_n\ebm
\]
as follows.  (Placing $\bbm$ between $\ba$ and $\bb$ is intended to suggest ``cutting $B$ by $\bbm$''.)
\blist
\item If there exists $i\in\{1,\dots,n\}$ such that $b_i<m_i$ (that is, if $\bbm \npreceq \bb$), then
\begin{equation}\label{gnomonformula1}
\bm \bb \\ \bbm \\ \ba \ebm = \bm \bb \\ \ba \ebm.
\end{equation}
\item If for all $i$, $m_i\leq b_i$ ($\bbm\preceq\bb$), let
\begin{equation}\label{defc}
c_i = \max\{a_i,m_i\}.
\end{equation}
Then
\begin{align}\label{gnomonformula2}
&\bm \bb \\ \bbm \\ \ba \ebm
= \bm m_1-1 & b_2&\cdots&b_n\\a_1&a_2&\cdots&a_n\ebm\\
&\qquad\sqcup \bm b_1&m_2-1&b_3&\cdots&b_n\\c_1&a_2&a_3&\cdots&a_n\ebm \notag\\
&\qquad\sqcup \bm b_1&b_2&m_3-1 &b_4&\cdots&b_n\\c_1&c_2&a_3&a_4&\cdots&a_n \ebm \notag\\
&\qquad\sqcup\cdots\sqcup \bm b_1&b_2&\cdots&b_{n-1}&m_n-1\\c_1&c_2&\cdots&c_n&a_n\ebm. \notag
\end{align}
\elist


The following lemma shows that the solution to the Block Subtraction Problem is a gnomon as defined above.

\blem\label{gnomonlem} For all nonnegative values of $a_i$, $b_i$, and $m_i$,
\begin{equation}\label{gnomonlem1}
\bm \bb \\ \ba \ebm \smallsetminus \langle \bbm\rangle
  = \bm \bb \\ \bbm \\ \ba \ebm.
\end{equation}
\elem

\bpf    Assume first that $b_i<m_i$ for some $i$. Then the ideal $\langle \bbm \rangle$ does not intersect the given block, and (\ref{gnomonlem1}) follows from (\ref{gnomonformula1}).

Next assume that $a_i<m_i\leq b_i$ for all $i$, so that $c_i=m_i$.  Then the original block can be decomposed into two subblocks by the hyperplane $\m_1=m_1$:
\begin{equation}\label{step1}
\bm \bb \\ \ba \ebm =
\bm m_1-1 & b_2&\cdots&b_n\\a_1&a_2&\cdots&a_n\ebm
\sqcup \bm b_1&b_2&b_3&\cdots&b_n\\c_1&a_2&a_3&\cdots&a_n\ebm.\notag
\end{equation}
(Note that $c_1$ occurs in the lower left corner of the second matrix where $m_1$ is expected, but these are equal.)  The first of these subblocks lies outside of $\langle \bx^{\bbm}\rangle$, and belongs to the difference we are computing; it equals the first block of (\ref{gnomonformula2}).  The second subblock is carried forward to the next step in a recursive process:  It is split by the hyperplane $\m_2=m_2$ into two subblocks, one having $a_2\leq \m_2<m_2$ and the other having $m_2\leq \m_2\leq b_2$.  Again, the first subblock is retained and the second carried forward to the third step, and so forth.

Finally, if $m_1\leq a_1$, then $m_1-1<a_1$ and the first subblock in (\ref{step1}) is empty.  Furthermore the lower left entry in the second matrix should be $a_1$ (rather than $m_1$, which would cause the block to include monomials that are not in the original block).  But in this case $c_1=a_1$, so (\ref{step1}) remains correct.  The same reasoning applies if the condition $m_i\leq a_i$ is encountered later in the recursion.  \epf

Now let $I=\langle \bbm^1,\dots,\bbm^r\rangle$ be a monomial ideal for which a block decomposition $B^1\sqcup\cdots\sqcup B^s$ is known, and let  $\langle \bbm \rangle$ be a principal monomial ideal.  Let $I'= I +
\langle \bbm\rangle$ be the sum of these ideals (in the usual sense), and note that $I'= \langle \bbm^1,\dots,\bbm^r,\bbm\rangle.$  We refer to this as {\bf adding a generator} to $I$.

\balg\label{addgenerator} A block decomposition for $S/I'$ is given by
\begin{equation}\label{addgenformula}
(B^1\smallsetminus \langle \bbm\rangle)\sqcup\cdots\sqcup(B^s\smallsetminus \langle \bbm\rangle).
\end{equation}
\ealg

\bpf  The standard monomials of $I'$ are the monomials that are standard for $I$ and, in addition, are not divisible by $\bx^{\bbm}$.  We remove the monomials divisible by $\bx^{\bbm}$ from each block $B^k$ of standard monomials for $I$ by finding the gnomons $B^k\smallsetminus \langle \bbm\rangle$.  Since these are disjoint, their union is a block decomposition.\epf

\brem\label{gnorder} For Section \ref{prime}, where the order of terms in a decomposition matters, we write the result of Algorithm \ref{addgenerator} in the form
\[
\bigsqcup_{k=1}^s \left( \bigsqcup_{i=1}^n B^{ki}\right)=
 (B^{11}\sqcup\dots\sqcup B^{1n})\sqcup\dots\sqcup (B^{s1}\sqcup\dots\sqcup B^{sn}).
\]
Here $B^{k1},\dots,B^{kn}$ are the terms of $B^k\smallsetminus \langle\bbm\rangle$ in the order shown in equation (\ref{gnomonformula2}), retaining any empty blocks, or, if (\ref{gnomonformula1}) applies, then $B^{k1}=B^k$ and $B^{ki}=\emptyset$ for $i=2,\dots,n$.  Retaining the empty sets allows a uniform notation in the proof of Theorem \ref{gnnice}.\erem

Let $I=\langle \bbm^1,\dots,\bbm^r \rangle$ be a monomial ideal given by its minimal generators.  The {\bf gnomon decomposition} of $S/I$, with respect to the standard order $x_1,\dots,x_n$ of the variables, is defined to be the result of applying Algorithm (\ref{gbd}) below.  The result may depend on the order in which the generators are indexed; this {\bf noncommutativity of the generators} is illustrated in the examples below.

\balg\label{gbd} Apply Algorithm \ref{addgenerator} repeatedly, beginning with $I=\emptyset$ and adding the generators $\bx^{\bbm^1},\dots,\bx^{\bbm^r}$ in that order.\ealg

Let $\pi$ be a permutation of the integers $(1,\dots,n)$, regarded as specifying an order $x_{\pi_1},\dots,x_{\pi_n}$ of the variables.  The $\pi$-{\bf gnomon decomposition} of $S/I$ is defined by changing variables temporarily to $y_i=x_{\pi_i}$, applying algorithm (\ref{gbd}) using the new variables, and then returning to the original variables.  (The generators are still taken in their indicated order.)

The blocks of the gnomon decomposition are unions of blocks of the elementary decomposition, because whenever a block is subdivided in Algorithm \ref{gbd}, the subdivision occurs along a portion of a hyperplane $\m_i=m^j_i$ for some $i$ and $j$.  The subdivisions for elementary blocks occur along these same hyperplanes (but along the entire hyperplane, not a portion of it.)

\bex\label{x1} For the example (\ref{firstexample0}), $I=\langle(3,9),(7,5)\rangle$.  According to Algorithm \ref{gbd}, we begin with the decomposition $\bm \infty&\infty\\0&0\ebm$ for $K[x,y]$ and add the generator $(3,9)$ using Algorithm \ref{addgenerator}, which in turn calls for using Lemma \ref{gnomonlem} to compute
\[
\bm \infty&\infty\\0&0\ebm \smallsetminus \langle (3,9)\rangle = \bm 2&\infty\\0&0\ebm \sqcup \bm \infty&8\\3&0\ebm.
\]
The first block contains the monomials that are not divisible by $x^3$, the second those that are divisible by $x^3$ but not by $y^9$.  Next we find
\[
\bm 2&\infty\\0&0\ebm \smallsetminus \langle (7,5)\rangle = \bm 2&\infty\\ 0&0\ebm
\]
since $2<7$ and $(7,5)$ is not in the block.  Finally
\[
\bm\infty&8\\3&0\ebm \smallsetminus \langle x^7y^5\rangle = \bm 6&8\\3&0\ebm\sqcup\bm \infty&4\\7&0\ebm.
\]
The result of these calculations is (\ref{firstexample2}).

The same problem works out as follows if we take $I=\langle (7,5),(3,9)\rangle$, that is, add the generators in the opposite order:
\begin{align*}
\bm \infty&\infty\\7&5\\0&0\ebm &= \bm 6&\infty\\0&0\ebm\sqcup\bm\infty&4\\7&0\ebm \\
\bm 6&\infty\\ 3&9\\ 0&0\ebm &= \bm 2&\infty\\0&0\ebm \sqcup \bm 6&8\\3&0\ebm\\
\bm \infty&4\\3&9\\7&0\ebm &= \bm\infty&4\\7&0\ebm.
\end{align*}
The last line is an instance of equation (\ref{gnomonformula1}).  The result is again (\ref{firstexample2}), so there is no noncommutativity in the generators in this example.  The $\pi$-gnomon decomposition with $\pi=(2,1)$ is (\ref{firstexample3}).\eex

\bex\label{x3} To show noncommutativity of generators, the gnomon decomposition of $K[x,y,z]/\langle (5,3,7),(10,6,2)\rangle$ is
\[
\bm 4&\infty&\infty\\0&0&0\ebm \sqcup \bm \infty&2&\infty\\5&0&0\ebm\sqcup \bm 9&\infty&6\\5&3&0\ebm \sqcup \bm\infty&5&6\\10&3&0\ebm \sqcup \bm\infty&\infty&1\\10&6&0\ebm,
\]
which is incompressible; the opposite order of generators gives instead
\begin{align*}
&\bm 4&\infty&\infty\\0&0&0\ebm \sqcup \bm 9&2&\infty\\5&0&0\ebm \sqcup \bm 9&\infty&6\\5&3&0\ebm\\
&\qquad \sqcup \bm \infty&2&\infty\\10&0&0\ebm \sqcup \bm \infty&5&6\\10&3&0\ebm \sqcup \bm\infty&\infty&1\\10&6&0\ebm.
\end{align*}
This is compressible (to the previous result), since
\[
\bm 9&2&\infty\\5&0&0\ebm \sqcup \bm \infty&2&\infty\\10&0&0\ebm = \bm \infty&2&\infty\\5&0&0\ebm.
\]
\eex

\bex\label{x4} The next example will be used in Remark \ref{x4cont} to illustrate Theorem \ref{gnnice}.  Let $S=K[x,y,z]$ and $I=\langle z^5,y^2z^3,x^3yz\rangle$.  Then
\begin{align*}
\bm \infty&\infty&\infty\\0&0&5\\0&0&0\ebm &= \bm \infty&\infty&4\\0&0&0\ebm\\
\bm \infty&\infty&4\\0&2&3\\0&0&0\ebm &= \bm \infty&1&4\\0&0&0\ebm\sqcup\bm\infty&\infty&2\\0&2&0\ebm\\
\bm \infty&1&4\\3&1&1\\0&0&0\ebm &= \bm 2&1&4\\0&0&0\ebm \sqcup \bm \infty&0&4\\3&0&0\ebm \sqcup \bm \infty&1&0\\3&1&0\ebm\\
\bm \infty&\infty&2\\3&1&1\\0&2&0\ebm &= \bm 2&\infty&2\\0&2&0\ebm \sqcup \quad\emptyset\quad \sqcup \bm \infty&\infty&0\\3&2&0\ebm.
\end{align*}
The gnomon decomposition for $S/I$ is the union of the last two lines.
\eex

We conclude this section with two lemmas that will be used in Section \ref{prime}.

\blem\label{intersect} With $c_i$ as in (\ref{defc}),
\begin{equation}\label{gnomonlem2}
\bm b_1&\cdots&b_n\\a_1&\cdots&a_n\ebm \cap \langle x_1^{m_1}\cdots x_n^{m_n}\rangle
  = \bm b_1&\cdots&b_n\\c_1&\cdots&c_n\ebm.
\end{equation}
\elem

\bpf If $m_i>b_i$ for any $i$, both sides are empty.  Otherwise, continue the reasoning in the proof of Lemma \ref{gnomonlem}.  Observe that the right-hand side of \ref{gnomonlem2} would be the next matrix added to (\ref{gnomonformula2}) if the pattern of (\ref{gnomonformula2}) is continued.\epf

\blem\label{gnlemma} Assume that $\bbm\preceq \bb$, so that (\ref{gnomonformula2}) is valid.  Let the blocks on the right-hand side of (\ref{gnomonformula2}) be labeled $B^1,\dots,B^n$.  Then for any $i=1,\dots,n$ we have
\[
B^i \sqcup\cdots\sqcup B^n \sqcup \bm \bb\\ \bc\ebm =
\bm  b_1 &\cdots &b_{i-1} &b_i &\cdots &b_n\\
     c_1 &\cdots &c_{i-1} &a_i &\cdots &a_n \ebm.
\]
\elem

\bpf Compute the unions successively starting from the right and working to the left, noticing that each union is between matrices that differ only in one column.  For instance, if $n=3$ the first step is
\[
\bm b_1&b_2&m_3-1\\c_1&c_2&a_3\ebm \sqcup \bm b_1&b_2&b_3\\c_1&c_2&c_3\ebm = \bm b_1&b_2&b_3\\c_1&c_2&a_3\ebm.
\]
If $c_3=m_3$ the result is clear; if $c_3=a_3$ then the left-hand block is empty and again the result is clear.  (The lemma and proof remain valid without assuming $\bbm\preceq\bb$, but in that case it does not apply to the gnomon decomposition, which is governed instead by (\ref{gnomonformula1}).)
\epf

\section{Subprime filtrations}\label{prime}

A {\bf filtration}  of $S$ by monomial ideals is a finite nested sequence
\begin{equation}\label{filt}
I=I_0\subset I_1\subset \cdots \subset I_s = S.
\end{equation}
This is often called a filtration of $S/I$, and it does induce an actual filtration of $S/I$ by subrings, namely $\{0\} \subset I_1/I\subset I_2/I\subset\cdots\subset S/I$.  If each set difference $\mon(I_k)\smallsetminus \mon(I_{k-1})$ is a block $B_k$, we call (\ref{filt}) a {\bf subprime filtration}.  If, further, each block $B_k$ is a Stanley block, we call (\ref{filt}) a {\bf prime filtration};  this is simpler than the usual definition in the literature, but is equivalent, as pointed out at the end of this section.  A block decomposition, {\it with the blocks indexed in a specific order}, will be called {\bf subprime} (or {\bf prime}) if it is associated in this way with a subprime (or prime) filtration.  In this section all block decompositions are taken as {\bf ordered block decompositions}.  This simplifies the discussion, in that we need not consider all possible orderings before deciding whether a given block decomposition is subprime.

The following equations hold trivially for subprime filtrations.  Quotients, direct sums, and isomorphisms are to be understood as quotients, direct sums, and isomorphisms of vector spaces, and $M$ is the set $\mon(S)\smallsetminus \mon(I)$ of standard monomials for $I$ in $S$.
\begin{align}
&I_k/I_{k-1} \cong \Span{B_k}\\
&S/I \cong \Span{B_1}\oplus\cdots\Span{B_s}\\
&M=B_1\sqcup\cdots\sqcup B_s\\
&\mon(I_{k-1})\sqcup B_k = \mon(I_k)\\
&I_{k-1}\oplus \Span(B_k)=I_k.
\end{align}
That is, with each subprime filtration of $S/I$ is associated a block decomposition in which the blocks are naturally ordered so that the $i$-th block $B_i=\mon(I_i)\smallsetminus \mon(I_{i-1})$, and when the block spaces are successively added (by direct sum) to the $I$ (considered as a vector space), the ideals of the filtration are reconstructed.  Thus the subprime filtration and the (ordered) block decomposition contain the same information.  Prime filtrations and (ordered) Stanley decompositions are related in the same way.

A convenient notation to exhibit the relationship between ideals and blocks in a subprime filtration is
\[
\begin{matrix} I=I_0 &\subset &I_1 &\subset &I_2 &\cdots &I_{s-1} &\subset &I_s=S.\\
{} &\vert &{} &\vert &{} &{} &{} &\vert \\
{} &B_1 &{} &B_2 &{} &{} &{} &B_s \end{matrix}
\]

For example, the block decomposition (\ref{bd1}) is associated with a subprime filtration as follows:
\[
\begin{matrix} \langle xy^3,x^3y\rangle &\subset &\langle xy \rangle &\subset &\langle y \rangle &\subset &\langle 1 \rangle,\\
{} &\vert &{} &\vert &{} &\vert &{}\\
{} &B_1 &{} &B_2 &{} &B_3 &{}
\end{matrix}
\]
with $B_1=B^3$, $B_2=B^2$, and $B_3=B^1$.

\brem\label{filtrationexample} Observe that the ideals in the filtration can be found by adding the inner corner of each block to the generators of the previous ideal.  For instance, adding the inner corner $(1,1)$ of $B_1=B^3$ to $\langle xy^3,x^3y\rangle$ gives $\langle xy^3,x^3y,xy\rangle=\langle xy\rangle$.  See Lemma \ref{obviouslemma} below.\erem

It is known that there exist Stanley decompositions that are not associated in this way with prime filtrations; thus there are also block decompositions that are not associated with subprime decompositions.  The standard example, introduced in \cite{McLagan-Smith}, is $S=K[x,y,z]$, $I=\langle xyz\rangle$,
\begin{equation}\label{MSexample}
S/I\cong K\oplus K[x,y]x\oplus K[y,z]y \oplus K[x,z]z.
\end{equation}
In this example, the Stanley blocks cannot be ordered in such a way that when added successively to $I$, they produce monomial ideals.  Here is our version of Soleyman-Jahan's theorem characterizing Stanley decompositions that come from prime filtrations.

\bthm\label{ASJ} An ordered block decomposition $B_1\sqcup\dots\sqcup B_s$ for $S/I$ is subprime if and only if $I\oplus \Span B_1 \oplus \cdots \oplus \Span(B_i)$ is a monomial ideal for each $i=1,\dots,s$.\ethm

\bpf The $\Rightarrow$ direction has already been proved.  For the $\Leftarrow$ direction, the hypothesis implies that $
\mon(I)\sqcup B_1\sqcup\cdots\sqcup B_k$ is a vector space basis for an ideal $I_k$ for each $k$, and that these ideals form a filtration.  It is equally clear that each $\mon(I_k)\smallsetminus \mon(I_{k-1})=B_k$, so the filtration is subprime.\epf

Note that when written as a block decomposition, the blocks of (\ref{MSexample}) are unions of blocks from the elementary decomposition of $S/I$, just as is the gnomon decomposition.  This shows that although the elementary decomposition itself is subprime (Theorem \ref{elemnice} below), the result of compressing an elementary decomposition is not in general subprime.  But the gnomon decomposition, which is a compression of the elementary decomposition, is subprime (Theorem \ref{gnnice}).

If $B$ is any block (\ref{defblock}), and $b_i<\infty$, we define the $i$th {\bf outer adjacent face} of $B$ to be the block
\[
F^i(B)=\bm b_1 &\dots &b_{i-1} &b_i +1 &b_{i+1} &\dots &b_n\\
           a_1 &\dots &a_{i-1} &b_i +1 &a_{i+1} &\dots &a_n\ebm.
\]
(This is not a subset of $B$, but is displaced by one step in the direction $i$ from the ordinary $i$th {\bf outer face} of $B$.)  If $b_i=\infty$ we set $F^i(B)=\emptyset$.

\blem\label{obviouslemma} Let $B$ be a block with inner corner $\ba$, and let $I$ be a monomial ideal disjoint from $B$.  Let $M=B\sqcup \mon(I)$.  Then $M$ spans a monomial ideal (or equivalently, $M= \mon(J)$ for some monomial ideal $J$) if and only if the outer adjacent faces of $B$ are contained in $I$.  In this case $B\sqcup \mon(I)=\mon(\langle \bx^{\ba}\rangle)\cup \mon(I)$ and $J=I+\langle \bx^{\ba}\rangle$.\elem

\bpf A set $M$ of monomials spans an ideal if and only if every monomial divisible by ($\succeq$) an element of $M$ belongs to $M$.  Let $M=B\sqcup \mon(I)$ and suppose that $M$ spans a monomial ideal.  Let $\bc\in F^i(B)$ for some $i$.  Then $\ba\prec\bc$, so $\bc\in M$.  But $\bc\notin B$, so $\bc\in I$.  Therefore $F^i(B)\subset I$ for each $i$.  Conversely, suppose that $F^i(B)\subset I$ for each $i$.  Since every element of $\langle\ba\rangle$ that is not in $B$ is divisible by some element of $F^i(B)$, and therefore belongs to $I$, we conclude that $\langle \ba\rangle \subset M$. Therefore $M= \mon(I) \cup \mon (\langle \ba\rangle)=\mon(I+\langle \bx^{\ba}\rangle)$, so $M$ spans an ideal.  \epf

We are now ready to prove that elementary and gnomon decompositions (with suitable orderings) are subprime.

\bthm\label{elemnice}  Let $I\subset S=K[x_1,\dots,x_n]$ be a monomial ideal and let $B^1\sqcup\dots\sqcup B^s$ be the elementary block decomposition of $S/I$, ordered according to Remark \ref{elemorder}.  Let $B_k=B^{s-k+1}$ for $k=1,\dots,s$; that is, reverse the ordering of the blocks.  Then $B_1\sqcup\dots\sqcup B_s$
is a subprime decomposition.\ethm

\bpf  We work in the superscript notation, restating the criterion of Theorem \ref{ASJ} as
\begin{equation}\label{subprimecriterion}
B^k \sqcup\dots\sqcup B^s\sqcup \mon(I)\quad \text{spans an ideal for each}\quad k=1,\dots,s.
\end{equation}
The proof begins with $k=s$ and works backward.  Let $A=\{\ba^1,\dots,\ba^s\}$ be the refined list of inner corners of elementary blocks for $S/I$, created in step \ref{refine} of Algorithm \ref{ebd}, in lexicographic order according to Remark \ref{elemorder}.  Since $\bbm \prec \bbm'$  implies $\bbm <_{lex} \bbm'$, $\ba^s$ is a maximal element of $A$ under $\prec$.

We claim that each outer adjacent face $F^i(B^s)$ is contained in $I$.  Let $\bc\in F^i(B^s)$.  Then $\ba^s\prec\bc$ and $\bc\notin B^s$.  But by Remark \ref{elemprop}, if $\bc\notin I$ then $\bc$ belongs to the block of $S/I$ whose inner corner is maximal among inner corners less than $\bc$.  But this inner corner must be $\ba^s$, so $\bc\in B^s$, a contradiction.  Therefore $\bc\in I$, proving the claim.

By Lemma \ref{obviouslemma} we conclude that
\[
B^s \sqcup I = I + \langle \ba^s \rangle = I'
\]
is an ideal.  Setting $A'=\{\ba^1,\dots,\ba^{s-1}\}$, we see that $A'$ is the list of inner corners of elementary blocks for $S/I'$, which has $B^1\sqcup\dots\sqcup B^{s-1}$ for its elementary decomposition.  Now the argument can be repeated, with $\ba^{s-1}$ as the maximal element.\epf

Lexicographic ordering is not the only ordering that gives a subprime decomposition.  It is only necessary to select a maximal element of $A$, remove it, and repeat.

\bthm\label{gnnice} The $\pi$-gnomon decomposition $B^1\sqcup\dots\sqcup B^s$ of $S/I$ is subprime, for any $\pi$, under the ordering $B_k=B^{s-k+1}$.\ethm

The proof of this theorem actually establishes more, namely, that if Algorithm \ref{addgenerator} is used to ``add a generator'' to a subprime block decomposition, the result will be subprime.

\bpf  Since any $\pi$-gnomon decomposition becomes a standard gnomon decomposition by renumbering the variables, it suffices to prove the standard case. The argument is by induction on the number $r$ of minimal generators of $I$.  If $r=1$, then
\[ \mon(I)=\bm \infty &\cdots&\infty\\ m_1&\cdots&m_n\ebm
\]
and the gnomon decomposition is
\[
\bm \infty&\cdots&\infty\\ m_1&\cdots&m_n\\0&\cdots&0\ebm.
\]
In this case Lemma \ref{gnlemma} shows that when the terms of the gnomon decomposition are added in reverse order to the ideal, an ideal is obtained at each stage.

Now make the induction hypothesis that the theorem is true for all monomial ideals $I$ having $r$ minimal generators.  Let $I'$ be a monomial ideal with $r+1$ generators, and write $I'=I + \langle \bbm\rangle$, where $I$ has $r$ generators.  Suppose that $B^1\sqcup\cdots\sqcup B^s$ is the gnomon decomposition for $S/I$.    The gnomon decomposition for $S/I'$ is
\[
(B^1\smallsetminus\langle\bbm\rangle)\sqcup\dots\sqcup(B^s\smallsetminus\langle\bbm\rangle)
 = \bigsqcup_{k=1}^s\left(\bigsqcup_{i=1}^n B^{ki}\right)
\]
in the notation of Remark \ref{gnorder}. We need to prove that for any $p\in\{1,\dots,s\}$ and $q\in\{1,\dots,n\}$,
\begin{equation}\label{toprove1}
(B^{pq}\sqcup\dots \sqcup B^{pn})\sqcup \bigsqcup_{k=p+1}^s\left(\bigsqcup_{i=1}^n B^{ki}\right) \sqcup \mon(I')
\end{equation}
spans an ideal.  From here on we assume $p$ and $q$ are fixed.  By the induction hypothesis,
\begin{equation}\label{indhyp}
B^k \sqcup\dots\sqcup B^s\sqcup \mon(I)
\end{equation}
spans an ideal for each $k=1,\dots,s$.  This will be used for $k=p+1$ and $k=p$.  For $k=p+1$, let $J$ be the ideal spanned by
\begin{equation}\label{defJ}
B^{p+1}\sqcup\dots\sqcup B^s\sqcup\mon(I).
\end{equation}
For $k=p$, (\ref{indhyp}) together with Lemma \ref{obviouslemma} implies that the outer adjacent faces of $B^p$ belong to $J$.

Let
\[
B^k=\bm \bb^k\\ \ba^k\ebm \qquad\text{and}\qquad c^k_i=\max\{a^k_i,m_i\}.
\]
By Lemma \ref{intersect} we have
\[
\bm \bb^k\\ \bc^k \ebm = B^k \cap \langle \bbm\rangle \subset I',
\]
so we can include these terms in (\ref{toprove1}) without changing the result, provided we change the last disjoint union to an ordinary union because at that point there is now an overlap.  Thus (\ref{toprove1}) equals
\begin{align}\label{toprove2}
&\left(B^{pq}\sqcup\dots \sqcup B^{pn}\sqcup \bm \bb^p\\ \bc^p\ebm \right)\\
&\qquad\qquad \sqcup \bigsqcup_{k=p+1}^s\left(\left(\bigsqcup_{i=1}^n B^{ki}\right)\sqcup \bm \bb^k\\ \bc^k\ebm \right) \cup \mon(I').\notag
\end{align}
By Lemmas \ref{gnomonlem} and \ref{intersect} we have
\[
\left(\bigsqcup_{i=1}^n B^{ki}\right)\sqcup \bm \bb^k\\ \bc^k\ebm = B^k.
\]
Since (\ref{defJ}) spans $J$ and $\mon(J)\cup\mon(I')=\mon(J+I')$,
\[
\left(\bigsqcup_{k=p+1}^s B^k\right) \cup \mon(I')
\]
spans $J'=J+I'$.  These calculations reduce (\ref{toprove2}) to
\begin{equation}\label{toprove3}
\left(B^{pq}\sqcup\dots \sqcup B^{pn}\sqcup \bm \bb^p\\ \bc^p\ebm \right) \cup \mon(J').
\end{equation}
The induction is completed by considering two cases.
\blist
\item If $\bbm^p\npreceq\bb^p$, then (\ref{gnomonformula1}) applies at the point when the generator $\bbm^p$ is added during the construction of the gnomon decomposition.  So by Remark \ref{gnorder}, $B^{p1}=B^p$ and $B^{pi}=\emptyset$ for $i=2,\dots,n$.  So (\ref{toprove3}) is trivial unless $p=1$.  But in that case, (\ref{toprove3}) reduces to $B^p\cup \mon(J')$.  The outer adjacent faces of $B^p$ belong to $J$, and therefore also to $J'$.  So (\ref{toprove3}) spans an ideal.

\item  If $\bbm^p\preceq\bb^p$, then Lemma \ref{gnlemma} implies
\begin{equation}\label{lastpiece}
B^{pq} \sqcup\cdots\sqcup B^{pn} \sqcup \bm \bb^p\\ \bc^p\ebm =
\bm  b^p_1 &\cdots &b^p_{q-1} &b^p_q &\cdots &b^p_n\\
     c^p_1 &\cdots &c^p_{q-1} &a^p_q &\cdots &a^p_n \ebm.
\end{equation}
The outer adjacent faces of this block are contained in those of
\[
B^p=\bm b^p_1&\cdots&b^p_n\\a^p_1&\cdots&a^p_n\ebm,
\]
but as noted above, these  belong to $J'$, so again, (\ref{toprove3}) spans an ideal.
\elist
\epf

\brem\label{x4cont} To find the subprime filtration associated with a gnomon decomposition, as in Remark \ref{filtrationexample}, one should delete any generators that are divisible by each new inner corner that is added.  The proof of Theorem \ref{gnnice} shows that when a sequence of blocks $B^{pn},\dots,B^{p1}$ is added  to $\mon(J')$, their inner corners decrease in $\prec$.  Therefore each new corner added will replace the previous one (and may also remove others).  For instance, in Example \ref{x4} the inner corners $(3,2,0),(0,2,0),(3,1,0),\newline (3,0,0),(0,0,0)$ are added to the starting ideal with generators $(0,0,5),\newline (0,2,3),(3,1,1)$.  The sequences $(3,2,0),(0,2,0)$ and $(3,1,0),(3,0,0),\newline (0,0,0)$ each arise from a single gnomon, so the inner corners decrease in $\prec$.   The subprime filtration is
\begin{align*}
&\langle (0,0,5),(0,2,3),(3,1,1)\rangle \subset \langle (0,0,5),(0,2,3),(3,1,1),(3,2,0)\rangle\\
&\qquad\subset\langle (0,0,5),(3,1,1),(0,2,0)\rangle
\subset \langle (0,0,5),(0,2,0), (3,1,0)\rangle\\
&\qquad \subset \langle (0,0,5),(0,2,0),(3,0,0)\rangle \subset \langle (0,0,0)\rangle.
\end{align*}
When $(0,2,0)$ is added, it eliminates both $(3,2,0)$, as expected, and also $(0,2,3)$, which is not predicted by this remark.
\erem

We conclude the section with some technical remarks concerning the usual definition of prime filtration.  Recall that a monomial ideal $P\subset S$ is prime if and only if it is generated by a subset $Y$ of the variables $x_1,\dots,x_n$.  Since $S$ is a $\Z^n$-graded algebra, and monomial ideals in $S$ inherit this multi-grading, quotients of monomial ideals are $\Z^n$-graded modules.  If $P$ is a monomial ideal and $\ba\in\N^n$, the symbol $S/P(-\ba)$ denotes $S/P$ with its $\Z^n$-grading shifted so that $(S/P)(-\ba)_{\bb} = (S/P)_{\bb+\ba}$.  If $P=\langle Y\rangle$ is a prime ideal, $S/P$ is isomorphic as a vector space to $K[X]$, where $X$ is the complementary set of variables to $Y$.   Then $S/P(-a)$ is isomorphic as a vector space to the Stanley space $K[X]\bx^{\ba}$.  Now (\ref{filt}) is defined to be a prime filtration if there exist monomial prime ideals $P_i$ and integer vectors $\ba_i\in\N^n$ such that $I_i/I_{i+1}$ is isomorphic as a $\Z^n$-graded $S$-module to $S/P(-\ba_i)$. It is therefore also isomorphic as a vector space to $K[X]\bx^{\ba}$, but this is not part of the definition.

Using these definitions, Soleyman-Jahan proved Theorem \ref{ASJ} above for the case of prime filtrations, with the blocks being Stanley blocks.  The difficulty of his proof seems to result from the fact that the association of $I_i/I_{i-1}$ with a particular $S/P_i(-\ba_i)$ is only specified by an isomorphism, and not by a concrete construction as we have done.  It seems to be the case that if $S/P(-\ba)$ is isomorphic as a $\Z^n$-graded $S$-module to $S/P'(-\ba')$, then $P=P'$ and $\ba=\ba'$.  In other words, $P_i$ and $\ba_i$ are uniquely determined by $I_i$ and $I_{i-1}$.  This observation (if correct) would bring his treatment closer to ours.  However, this uniqueness is only possible if the isomorphism $I_i/I_{i-1}\cong S/P_i(-\ba_i)$ is interpreted as an $S$-module homomorphism; neither algebra homomorphism nor vector space homomorphism is strong enough to achieve the desired uniqueness.  With our approach there is no need to bring in the $\Z^n$ grading or the $S$-module structure.

\bibliography{Murdock}

\end{document}